\documentclass[12pt]{amsart}

\usepackage{array}

\usepackage{pdflscape}
\usepackage{graphicx}
\usepackage{float}
\usepackage{caption}

\usepackage{amssymb,amsmath}
\usepackage{amssymb,latexsym}
\usepackage{amsfonts}
\usepackage{amssymb}
\usepackage{stmaryrd}
\usepackage{longtable}
\usepackage{amsmath, geometry, amssymb} 
\usepackage{chngcntr}
\numberwithin{equation}{section}
\newtheorem{theorem}{Theorem}[section]

\newcommand{\sqr}[2]{{\vcenter{\vbox{\hrule height#2pt
                \hbox{\vrule width#2pt height#1pt \kern#1pt
                \vrule width#2pt}\hrule height#2pt}}}}

\newcommand{\mb}{\mbox}

\newcommand{\beq}{\begin{equation}}
\newcommand{\eeq}{\end{equation}}
\newcommand{\beqar}{\begin{eqnarray}}
\newcommand{\eeqar}{\end{eqnarray}}
\def\beqars{\begin{eqnarray*}}
\def\eeqars{\end{eqnarray*}}
\def\eop{\hfill\mb{$\hspace{12pt}\vrule height 7pt width 6pt depth 0pt$}}

\newcommand{\pmd}{\hspace{-3mm} \pmod}

\newcommand{\smod}[1]{\hspace{-1mm} \pmod{#1}}

\newcommand{\qu}[2]{\Bigl({\frac{#1}{#2}}\Bigr) }
\newcommand{\dqu}[2]{\ds{\qu{#1}{#2}}}

\def \ds{\displaystyle}

\newcommand{\nn}{\mathbb{N}}
\newcommand{\zz}{\mathbb{Z}}
\newcommand{\qq}{\mathbb{Q}}
\newcommand{\cc}{\mathbb{C}}

\newcommand{\hh}{\mathbb{H}}

\allowdisplaybreaks

\begin{document}

\title{Representations by octonary quadratic forms\\
 with coefficients $1$, $2$, $3$ or $6$}  

\author{Ay\c{s}e Alaca and  M. Nesibe Kesicio\u{g}lu}

\maketitle

\markboth{AY\c{S}E ALACA AND  M. NESIBE KESICIO\u{G}LU }
{REPRESENTATIONS BY OCTONARY QUADRATIC FORMS}

\begin{abstract}

Using modular forms we determine formulas for the number of representations of a positive integer
by diagonal octonary quadratic forms  with coefficients $1$, $2$, $3$ or $6$.\\


\end{abstract}

\section{Introduction}

Let $\nn$, $\nn_0$, $\zz$, $\qq$ and $\cc$ denote the sets of
positive integers, nonnegative integers, integers, rational numbers and complex numbers respectively.
For $a_1,\ldots ,a_8 \in \nn$  and $n\in \nn_0$, we define
\beqars
N(a_1,\ldots , a_8 ;n):=
{\rm card}\{(x_1,\ldots, x_8)\in \zz^8 \mid n= a_1x_1^2 +\cdots + a_8 x_8^2 \} . \hspace{-15mm}
\eeqars
Clearly $N(a_1,\ldots , a_8 ;0) =1$. Without loss of generality we may suppose that
\beqars
a_1 \leq \cdots \leq a_8 \mbox{~~and~~} {\rm gcd(}a_1, \ldots ,a_8 {\rm )}=1.
\eeqars
Formulae for $N(a_1,\ldots , a_8 ;n)$ for the octonary quadratic forms
\beqar
&&\hspace{-10mm}(x_1^2 + \cdots + x_{i}^2) + 2(x_{i+1}^2 + \cdots + x_{i+j}^2) \nonumber \\
&&\hspace{10mm} +  3(x_{i+j+1}^2 + \cdots +  x_{i+j+k}^2 )
+  6(x_{i+j+k+1}^2 + \cdots  + x_{i+j+k+l}^2)
\eeqar
with $i+j+k+l =8$
under  the conditions
\beqars
 i \equiv j \equiv k \equiv l \equiv 0 \pmd 2 \mbox{~~or~~} i \equiv j \equiv k \equiv l \equiv 1 \pmd 2
\eeqars
appeared in literature. See \cite{AAW4}, \cite{AK}, \cite{AK2}, \cite{SW},  \cite{jacobi}, \cite{kokluce-1} and \cite{ramak}.
For convenience, we write $(i,j,k,l)$ to denote an octonary quadratic form given by (1.1), and
we write $N(1^{i}, 2^{j}, 3^{k}, 6^{l} ;n)$ to denote the number of representations of $n$ by the octonary
quadratic form $(i,j,k,l)$.

In this paper we determine a formula for $N(1^{i}, 2^{j}, 3^{k}, 6^{l} ;n)$ for each of the octonary quadratic forms $(i,j,k,l)$ when
some of the $i,j,k$ or $l$ have different parities from the others. There are one hundred and twelve such cases, and all of them are listed in Table 1. 
This paper completes the representations of a positive integer by diagonal octonary quadratic forms  with coefficients $1$, $2$, $3$ or $6$.

\begin{table}[H]
\captionof{table}{Octonary quadratic forms $(i,j,k,l)$}
\begin{center}
\begin{tabular}{|c|c|c|} \hline
 $(i,j,k,l)$ & $(i,j,k,l)$ & $(i,j,k,l)$ \\ \hline
 $(0,2,1,5)$, $(2,0,5,1)$  & $(0,1,2,5)$, $(2,1,4,1)$   & $(0,1,1,6)$, $(2,1,1,4)$   \\
 $(0,2,3,3)$, $(2,2,1,3)$  & $(0,1,4,3)$, $(2,3,0,3)$   & $(0,1,3,4)$, $(2,1,3,2)$     \\
 $(0,2,5,1)$, $(2,2,3,1)$  & $(0,1,6,1)$, $(2,3,2,1)$   & $(0,1,5,2)$, $(2,1,5,0)$     \\
  $(0,4,1,3)$, $(2,4,1,1)$  & $(0,3,2,3)$, $(2,5,0,1)$   & $(0,1,7,0)$, $(2,3,1,2)$      \\
 $(0,4,3,1)$, $(3,1,0,4)$  & $(0,3,4,1)$, $(3,0,1,4)$   & $(0,3,1,4)$, $(2,3,3,0)$     \\
 $(0,6,1,1)$, $(3,1,2,2)$  & $(0,5,2,1)$, $(3,0,3,2)$   & $(0,3,3,2)$, $(2,5,1,0)$    \\
 $(1,1,0,6)$, $(3,1,4,0)$  & $(1,0,1,6)$, $(3,0,5,0)$   & $(0,3,5,0)$, $(3,0,0,5)$    \\
 $(1,1,2,4)$, $(3,3,0,2)$  & $(1,0,3,4)$, $(3,2,1,2)$   & $(0,5,1,2)$, $(3,0,2,3)$   \\
 $(1,1,4,2)$, $(3,3,2,0)$  & $(1,0,5,2)$, $(3,2,3,0)$   & $(0,5,3,0)$, $(3,0,4,1)$    \\
 $(1,1,6,0)$, $(3,5,0,0)$ & $(1,0,7,0)$, $(3,4,1,0)$   & $(0,7,1,0)$, $(3,2,0,3)$   \\
 $(1,3,0,4)$, $(4,0,1,3)$  & $(1,2,1,4)$, $(4,1,0,3)$   & $(1,0,0,7)$, $(3,2,2,1)$    \\
 $(1,3,4,0)$, $(4,0,3,1)$  & $(1,2,3,2)$, $(4,1,2,1)$   & $(1,0,2,5)$, $(3,4,0,1)$    \\
 $(1,3,2,2)$, $(4,2,1,1)$  & $(1,2,5,0)$, $(4,3,0,1)$   & $(1,0,4,3)$, $(4,1,1,2)$     \\
 $(1,5,0,2)$, $(5,1,0,2)$   & $(1,4,1,2)$, $(5,0,1,2)$   & $(1,0,6,1)$, $(4,1,3,0)$    \\
 $(1,5,2,0)$,  $(5,1,2,0)$   & $(1,4,3,0)$, $(5,0,3,0)$   & $(1,2,0,5)$, $(4,3,1,0)$ \\
 $(1,7,0,0)$,  $(5,3,0,0)$   & $(1,6,1,0)$, $(5,2,1,0)$   & $(1,2,2,3)$, $(5,0,0,3)$ \\
 $(2,0,1,5)$,  $(6,0,1,1)$   & $(2,1,0,5)$, $(6,1,0,1)$   & $(1,2,4,1)$, $(5,0,2,1)$ \\
 $(2,0,3,3)$, $(7,1,0,0)$   & $(2,1,2,3)$, $(7,0,1,0)$   & $(1,4,0,3)$, $(5,2,0,1)$\\
                                      &                                       & $(1,4,2,1)$, $(6,1,1,0)$    \\
                                      &                                       & $(1,6,0,1)$, $(7,0,0,1)$ \\
\hline
\end{tabular}
\end{center}
\end{table}

\section{Preliminary results}

For  $q\in \cc$ with $|q|<1$
we define
\beqar
 F(q):=  \prod_{n=1}^{\infty} (1-q^{n}).
\eeqar
Ramanujan's theta function  $\varphi (q)$ is defined by
\beqars
\ds \varphi (q) = \sum_{n=-\infty}^{\infty} q^{n^2}.
\eeqars
We note that
\beqar
\sum_{n=0}^{\infty} N(a_1,\ldots , a_8 ;n) q^n = \varphi (q^{a_1}) \cdots \varphi (q^{a_8}) .
\eeqar
The infinite product representation of $\varphi (q)$ is due to Jacobi \cite{jacobi}, namely
\beqar
\varphi (q) = \frac{F^5 (q^2)}{F^2 (q) F^2 (q^4)} .
\eeqar
The Dedekind eta function $\eta (z)$ is defined on the upper half plane $\hh = \{ z \in \cc \mid \mbox{\rm Im}(z) >0 \}$
by the product formula
\beqar
\eta (z) = e^{\pi i z/12} \prod_{n=1}^{\infty} (1-e^{2\pi inz}).
\eeqar
Throughout the remainder of the paper we take
$q=q(z):=e^{2\pi i z} ~\mbox{  with } z\in \hh$
so that  $|q| <1$. By (2.1) and (2.4) we have
\beqar
\eta (z) = q^{1/24} \prod_{n=1}^{\infty} (1-q^{n}) = q^{1/24} F(q).
\eeqar
An eta quotient is defined to be a finite product of the form
\beqar
f(z) = \prod_{\delta } \eta^{r_{\delta}} ( \delta z),
\eeqar
where $\delta$ runs through a finite set of positive integers and the exponents $r_{\delta}$ are non-zero integers.
By taking $N$ to be the least common multiple of the $\delta$'s we can write the eta quotient {\rm (2.6)}  as
\beqar
f(z) = \prod_{\delta \mid N} \eta^{r_{\delta}} ( \delta z) ,
\eeqar
where some of the exponents $r_{\delta}$ may be $0$.

Let $N\in \nn$ and $\chi$ be a Dirichlet character of modulus dividing $N$ and $\Gamma_0(N)$ the modular subgroup defined by
\beqars
\Gamma_0(N) = \left\{ \left(
\begin{array}{cc}
a & b \\
c & d
\end{array}
\right) \Big | \;  a,b,c,d\in \zz ,~ ad-bc = 1,~c \equiv 0 \pmd {N}
\right\} .
\eeqars
Let $k\in \zz$. We write $M_k(\Gamma_0(N),\chi)$ to denote the space of modular forms of weight $k$ with multiplier system $\chi$ for $\Gamma_0(N)$, and $E_{k}(\Gamma_0(N),\chi)$ and $S_{k}(\Gamma_0(N),\chi)$ to denote the subspaces of Eisenstein forms and cusp forms of $M_{k}(\Gamma_0(N),\chi)$, respectively. It is known (see for example \cite[p. 83]{stein} and \cite{Serre}) that
\beqar
M_{k}(\Gamma_0(N),\chi)=E_{k}(\Gamma_0(N),\chi)\oplus S_{k}(\Gamma_0(N),\chi).
\eeqar
We use the following theorem to determine if certain eta quotients are modular forms.
See \cite[Theorem 5.7, p. 99]{Kilford}, \cite[Corollary 2.3, p. 37]{Kohler},
\cite[p. 174]{GordonSinor} and \cite{Ligozat}.

\begin{theorem}{\bf\em{(Ligozat)}}
Let $N\in \nn$. Let
$\displaystyle  f(z) = \prod_{1 \leq \delta \mid N} \eta^{r_{\delta}} ( \delta z)$ be an eta quotient and
$s=\displaystyle \prod_{1\leq \delta \mid N} \delta^{r_\delta}$. Suppose that
$k=\displaystyle \frac{1}{2} \sum_{1\leq \delta \mid N}r_\delta $  is an integer.
If $f(z)$ satisfies the conditions

{\em (L1)~~} $\displaystyle \sum_{ 1\leq  \delta \mid N} \delta \cdot r_{\delta} \equiv 0 \smod {24}$,\\

{\em (L2)~} $\displaystyle \sum_{ 1 \leq \delta \mid N} \frac{N}{\delta} \cdot r_{\delta} \equiv 0 \smod {24}$, \\

{\em (L3)~} for each $d \mid N$, $\displaystyle \sum_{1 \leq \delta \mid N} \frac{ \gcd (d, \delta)^2 \cdot r_{\delta} }{\delta} \geq 0 $, \\

\noindent
then $f(z)\in  M_k (\Gamma_0(N),\chi)$, where $\chi$ is given by
\begin{eqnarray}
 \chi(m)=\Big(\frac{(-1)^k s}{m} \Big )      \nonumber
\end{eqnarray}

In addition to the above conditions if $f(z)$ also satisfies the condition

{\em (L4)~} for each $d \mid N$, $\displaystyle \sum_{1 \leq \delta \mid N} \frac{ \gcd (d, \delta)^2 \cdot r_{\delta} }{\delta} > 0 $, \\

\noindent
then $f(z) \in S_k(\Gamma_0(N), \chi)$.
\end{theorem}
\vspace{3mm}

Let $\psi_1$ and $\psi_2$ be Dirichlet characters. For $n\in \nn$ we define $\sigma_{(3,\psi_1,\psi_2)}(n)$ by

\beqar
\sigma_{(3,\psi_1,\psi_2)}(n):= \sum_{1\leq m \mid n} \psi_2(m)\psi_1(n/m)m^3.
\eeqar
If $n\not\in \nn$, we set $\sigma_{(3,\psi_1,\psi_2)}(n)=0$. If $\psi_1$ and $\psi_2$ are trivial characters then
$\sigma_{(3,\psi_1,\psi_2)}(n)$ coincides with the sum of divisors function
\beqars
\sigma_3(n) = \sum_{1\leq m \mid n} m^3.
\eeqars
Let $\chi_0$ denote the trivial character. For $m\in \zz$, we define six characters by
\beqar
\left\{\begin{array}{l}
\chi_1 ( m)=\dqu{-8}{m},~\chi_2 (m) =\dqu{-4}{m},~\chi_3 (m) =\dqu{-3}{m},\\[3mm]
\chi_4 (m) =\dqu{8}{m},~\chi_5 (m) =\dqu{12}{m},~\chi_6 (m) =\dqu{24}{m}.
\end{array}\right.
\eeqar
We define the Eisenstein series $E_{4,\chi_0,\chi_4} (q)$, $E_{4,\chi_0,\chi_5} (q)$, $E_{4,\chi_0,\chi_6} (q)$, $E_{4,\chi_1,\chi_3} (q)$,
$E_{4,\chi_2,\chi_3} (q)$, $E_{4,\chi_3,\chi_1} (q)$, $E_{4,\chi_4,\chi_0} (q)$, $E_{4,\chi_5,\chi_0} (q)$ and , $E_{4,\chi_6,\chi_0} (q)$   by
\beqar
&& E_{4,\chi_0,\chi_4} (q):=\displaystyle\frac{11}{2}+\sum_{n=1}^{\infty} \sigma_{(3,\chi_0,\chi_4)}(n)q^n,\\
&& E_{4,\chi_0,\chi_5} (q):=\displaystyle23+\sum_{n=1}^{\infty} \sigma_{(3,\chi_0,\chi_5)}(n)q^n,\\
&& E_{4,\chi_0,\chi_6} (q):=\displaystyle261+\sum_{n=1}^{\infty} \sigma_{(3,\chi_0,\chi_6)}(n)q^n,\\
&& E_{4,\chi_1,\chi_3} (q):=\displaystyle\sum_{n=1}^{\infty} \sigma_{(3,\chi_1,\chi_3)}(n)q^n,\\
&& E_{4,\chi_2,\chi_3} (q):=\displaystyle\sum_{n=1}^{\infty} \sigma_{(3,\chi_2,\chi_3)}(n)q^n,\\
&& E_{4,\chi_3,\chi_1} (q):=\displaystyle\sum_{n=1}^{\infty} \sigma_{(3,\chi_3,\chi_1)}(n)q^n,\\
&& E_{4,\chi_3,\chi_2} (q):=\displaystyle\sum_{n=1}^{\infty} \sigma_{(3,\chi_3,\chi_2)}(n)q^n,\\
&& E_{4,\chi_4,\chi_0} (q):=\displaystyle\sum_{n=1}^{\infty} \sigma_{(3,\chi_4,\chi_0)}(n)q^n,\\
&& E_{4,\chi_5,\chi_0} (q):=\displaystyle\sum_{n=1}^{\infty} \sigma_{(3,\chi_5,\chi_0)}(n)q^n,\\
&& E_{4,\chi_6,\chi_0} (q):=\displaystyle\sum_{n=1}^{\infty} \sigma_{(3,\chi_6,\chi_0)}(n)q^n.
\eeqar


\section{Main results}

We define the eta quotients  $A_k (q)$, $B_k(q)$ and $C_k(q)$ as in Table 2. 
\begin{table}
\setlength\extrarowheight{4pt}
\caption{Eta quotients $A_k (q), B_k(q)$ and $C_k(q)$}
\begin{tabular}{|c|c|c|c|} \hline
$k$  & $A_k(q)$ & $B_k(q)$ & $C_k(q)$ \\ \hline
$1$  & $\frac{\eta^{2}(3z)\eta^{2}(4z)\eta^{5}(6z)\eta^{2}(8z)}{\eta^{3}(12z)}$                      & $\frac{\eta^{4}(3z)\eta^{2}(6z)\eta^{3}(8z)}{\eta(24z)}$                             & $\frac{\eta^{7}(6z)\eta^{3}(8z)\eta^{3}(12z)}{\eta^{2}(3z)\eta^{3}(24z)}$            \\[1mm]
$2$  & $\frac{\eta^{2}(3z)\eta^{5}(4z)\eta(6z)\eta(24z)}{\eta(8z)}$                                            & $\frac{\eta^{4}(6z)\eta^{2}(8z)\eta^{5}(12z)}{\eta(4z)\eta^{2}(24z)}$                & $\frac{\eta^{2}(3z)\eta^{7}(4z)\eta^{4}(12z)}{\eta^{3}(6z)\eta^{2}(8z)}$      \\[1mm]
$3$  & $\frac{\eta^{5}(2z)\eta(3z)\eta(12z)\eta^{2}(24z)}{\eta(z)}$                                            & $\frac{\eta^{4}(3z)\eta^{2}(4z)\eta^{2}(6z)\eta^{3}(24z)}{\eta(8z)\eta^{2}(12z)}$    & $\frac{\eta^{2}(3z)\eta(6z)\eta^{6}(8z)\eta^{2}(12z)}{\eta^{3}(4z)}$        \\[1mm]
$4$  & $\frac{\eta^{2}(3z)\eta(6z)\eta^{4}(8z)\eta(12z)\eta^{2}(24z)}{\eta^{2}(4z)}$                & $\frac{\eta^{4}(6z)\eta^{3}(8z)\eta^{3}(24z)}{\eta^{2}(12z)}$                & $\frac{\eta^{2}(3z)\eta^{3}(8z)\eta^{5}(12z)\eta(24z)}{\eta^{3}(6z)}$                \\[1mm]
$5$  & $\frac{\eta^{2}(3z)\eta(4z)\eta(8z)\eta^{4}(12z)\eta^{3}(24z)}{\eta^{3}(6z)}$                & $\eta^{3}(4z)\eta(12z)\eta^{4}(24z)$                                                 & $\frac{\eta^{2}(3z)\eta(6z)\eta^{2}(8z)\eta^{4}(24z)}{\eta(4z)}$                     \\[1mm]
$6$  & $\frac{\eta^{2}(3z)\eta(6z)\eta^{6}(24z)}{\eta(12z)}$                                                   & $\frac{\eta^{2}(4z)\eta^{4}(6z)\eta^{7}(24z)}{\eta(8z)\eta^{4}(12z)}$                & $\frac{\eta^{2}(3z)\eta^{2}(4z)\eta^{3}(12z)\eta^{5}(24z)}{\eta^{3}(6z)\eta(8z)}$    \\[1mm]
$7$  & $\frac{\eta^{2}(3z)\eta^{3}(4z)\eta^{2}(12z)\eta^{7}(24z)}{\eta^{3}(6z)\eta^{3}(8z)}$                   & $\frac{\eta^{5}(4z)\eta^{8}(24z)}{\eta^{4}(8z)\eta(12z)}$                            & $\frac{\eta^{2}(3z)\eta(4z)\eta(6z)\eta^{8}(24z)}{\eta^{2}(8z)\eta^{2}(12z)}$        \\[1mm]
$8$  & $\frac{\eta^{5}(2z)\eta(3z)\eta^{3}(12z)\eta^{8}(24z)}{\eta(z)\eta^{2}(4z)\eta^{4}(6z)\eta^{2}(8z)}$    & $\frac{\eta^{4}(2z)\eta^{7}(24z)}{\eta(8z)\eta^{2}(12z)}$                            & $\frac{\eta(z)\eta(6z)\eta(12z)\eta^{8}(24z)}{\eta(3z)\eta^{2}(8z)}$                 \\[1mm]
$9$  & $\frac{\eta^{2}(3z)\eta^{6}(4z)\eta^{5}(12z)}{\eta^{3}(6z)\eta^{2}(24z)}$                               & $$                                                                                   & $\frac{\eta^{2}(2z)\eta^{6}(3z)\eta(4z)\eta^{2}(8z)}{\eta^{3}(6z)}$                  \\[1mm]
$10$ & $\frac{\eta^{2}(4z)\eta^{7}(6z)\eta^{2}(8z)\eta^{4}(24z)}{\eta^{2}(3z)\eta^{5}(12z)}$          & $$                                                                                   & $\frac{\eta^{2}(3z)\eta^{3}(4z)\eta^{5}(6z)\eta^{2}(24z)}{\eta^{4}(12z)}$            \\[1mm]
\hline
\end{tabular}
\end{table}
\noindent
We define the integers $a_k(n), b_k(n)$  and $c_k(n)$ ($n\in\nn$) by
\beqar
&&A_k(q)=\sum_{n=1}^{\infty} a_k(n) q^n,~~(1\leq k \leq 10),\\
&&B_k(q)=\sum_{n=1}^{\infty} b_k(n) q^n,~~1\leq k \leq 8),\\
&&C_k(q)=\sum_{n=1}^{\infty} c_k(n) q^n,~~(1\leq k \leq 10).
\eeqar
We note that the eta quotients in Table 2 are constructed by using MAPLE 
in a way that they satisfy the conditions of Theorem 2.1 for $N=24$ and $k=4$. 
We deduce from \cite[Sec. 6.3, p. 98]{stein} that
\begin{eqnarray}
&&\dim(S_4(\Gamma_0(24),\chi_4)=10, \\
&&\dim(S_4(\Gamma_0(24),\chi_5)=8,\\
&&\dim(S_4(\Gamma_0(24),\chi_6)=10.
\end{eqnarray}

\begin{theorem} 
{\em (a)~} $\{A_{k}(q)\mid 1\leq k\leq 10\}$ is a basis for  $S_4(\Gamma_0(24),\chi_4)$.

{\em (b)~}
$\{E_{4,\chi_0,\chi_4}(q^t), E_{4,\chi_4,\chi_0}(q^t) \mid t=1,3\}$
is a basis for $E_4(\Gamma_0(24),\chi_4)$.

{\em (c)~}
$\{E_{4,\chi_0,\chi_4}(q^t), E_{4,\chi_4,\chi_0}(q^t) \mid t=1,3\}\cup\{A_{k}(q)\mid 1\leq k\leq 10\}$ is a basis for $M_4(\Gamma_0(24),\chi_4)$.
\end{theorem}

{\bf Proof.}
{(a)} By Theorem 2.1,  $A_{k}(q)$  ($1\leq k\leq 10$) are in $S_4(\Gamma_0(24),\chi_4)$. Then the assertion follows from (3.4) as  there  is no linear relationship among the $A_k(q)$ $(1\leq k \leq 10)$.

{(b)} The assertion follows from  \cite[Theorem 5.9, p. 88]{stein} with $\epsilon=\chi_4$ and $\chi, \psi \in \{\chi_0, \chi_4\}$.

{(c)} The assertion follows from (a), (b) and (2.8).
\eop

\begin{theorem} 
{\em (a)} $\{B_{k}(q)\mid 1\leq k\leq 8\}$ is a basis for  $S_4(\Gamma_0(24),\chi_5)$.

{\em (b)}
$\{E_{4,\chi_0,\chi_5}(q^t), E_{4,\chi_2,\chi_3}(q^t), E_{4,\chi_3,\chi_2}(q^t), E_{4,\chi_5,\chi_0}(q^t) \mid t=1,2\}$
is a basis for $E_4(\Gamma_0(24),\chi_5)$.

{\em (c)}
$\{E_{4,\chi_0,\chi_5}(q^t), E_{4,\chi_2,\chi_3}(q^t), E_{4,\chi_3,\chi_2}(q^t), E_{4,\chi_5,\chi_0}(q^t) \mid t=1,2\}\cup \{B_{k}(q)\mid 1\leq k\leq 8\}$
is a basis for $M_4(\Gamma_0(24),\chi_5)$.
\end{theorem}

{\bf Proof.}
{(a)} By Theorem 2.1,  $B_{k}(q)$  ($1\leq k\leq 8$) are in $S_4(\Gamma_0(24),\chi_5)$. Then the assetion follows from (3.5) as there is no linear relationship among the $B_k(q)$ $(1\leq k \leq 8)$.

{(b)} The assertion follows from  \cite[Theorem 5.9, p. 88]{stein} with $\epsilon=\chi_5$ and $\chi, \psi \in \{\chi_0, \chi_2, \chi_3, \chi_5\}$.

{(c)} The assertion follows from (a), (b) and (2.8).
\eop

\begin{theorem} 
{\em (a)} $\{C_{k}(q)\mid 1\leq k\leq 10\}$ is a basis for  $S_4(\Gamma_0(24),\chi_6)$.

{\em (b)}
$\{E_{4,\chi_0,\chi_6}(q), E_{4,\chi_1,\chi_3}(q), E_{4,\chi_3,\chi_1}(q), E_{4,\chi_6,\chi_0}(q)\}$
is a basis for $E_4(\Gamma_0(24),\chi_6)$.

{\em (c)}
$\{E_{4,\chi_0,\chi_6}(q), E_{4,\chi_1,\chi_3}(q), E_{4,\chi_3,\chi_1}(q), E_{4,\chi_6,\chi_0}(q)\}\cup\{C_{k}(q)\mid 1\leq k\leq 10\}$ is a basis for $M_4(\Gamma_0(24),\chi_6)$.
\end{theorem}

{\bf Proof.}
{(a)} By Theorem 2.1,  $C_{k}(q)$  ($1\leq k\leq 8$) are in $S_4(\Gamma_0(24),\chi_6)$. Then the assetion follows from (3.6) as there is no linear relationship among the $C_k(q)$ $(1\leq k \leq 10)$.

{(b)} The assertion follows from  \cite[Theorem 5.9, p. 88]{stein} with $\epsilon=\chi_6$ and $\chi, \psi \in \{\chi_0, \chi_1, \chi_3, \chi_6\}$.

{(c)} The assertion follows from (a), (b) and (2.8).
\eop

\begin{theorem} 
Let  $(i,j,k,l)$ be any of the octonary quadratic forms listed in the first column of {\em Table 1}. Then
\beqars
\varphi^{i}(q)\varphi^{j}(q^2)\varphi^{k}(q^3)\varphi^{l}(q^6) &=& x_1E_{4,\chi_0,\chi_4}(q) + x_2E_{4,\chi_4,\chi_0}(q) +  x_3E_{4,\chi_0,\chi_4}(q^3) + x_4E_{4,\chi_4,\chi_0}(q^3) \\
&&+ \sum_{i=1}^{10}y_iA_i(q),
\eeqars
where $x_i$  $(1\leq i \leq 4)$ and $y_i$ $(1\leq i \leq 10)$ are listed in {\em Table 3}.
\end{theorem}

{\bf Proof.} By (2.3), (2.5) and Theorem 2.1, we have $\varphi^{i}(q)\varphi^{j}(q^2)\varphi^{k}(q^3)\varphi^{l}(q^6) \in M_4(\Gamma_0(24),\chi_4)$.
Appealing to Theorem 3.1(c) and using MAPLE,  we obtain the asserted results.
\eop

\begin{theorem} 
Let  $(i,j,k,l)$ be any of the octonary quadratic forms  listed in the second column of Table 1. Then
\beqars
\varphi^{i}(q)\varphi^{j}(q^2)\varphi^{k}(q^3)\varphi^{l}(q^6)
&=& x_1E_{4,\chi_0,\chi_5}(q) + x_2E_{4,\chi_2,\chi_3}(q) + x_3E_{4,\chi_3,\chi_2}(q) +  x_4E_{4,\chi_5,\chi_0}(q)   \\
&&+ x_5E_{4,\chi_0,\chi_5}(q^2) + x_6E_{4,\chi_3,\chi_2}(q^2) + \sum_{i=1}^{8}y_iB_i(q),
\eeqars
where $x_i$ $(1\leq i \leq 6)$ and $y_i$ $(1\leq i \leq 8)$ are listed in {\em Table 4}.
\end{theorem}

{\bf Proof.} By (2.3), (2.5) and Theorem 2.1, we have $\varphi^{i}(q)\varphi^{j}(q^2)\varphi^{k}(q^3)\varphi^{l}(q^6) \in M_4(\Gamma_0(24),\chi_5)$.
Appealing to Theorem 3.2(c) and using MAPLE,  we obtain the asserted results.
\eop

\begin{theorem} 
Let  $(i,j,k,l)$ be any of the octonary quadratic forms  listed in the third column of Table 1. Then
\beqars
\varphi^{i}(q)\varphi^{j}(q^2)\varphi^{k}(q^3)\varphi^{l}(q^6)& =& x_1E_{4,\chi_0,\chi_6}(q) + x_2E_{4,\chi_1,\chi_3}(q) + x_3E_{4,\chi_3,\chi_1}(q) +  x_4E_{4,\chi_6,\chi_0}(q)\\
&& + \sum_{i=1}^{10}y_iC_i(q),
\eeqars
where $x_i$ $(1\leq i \leq 4)$ and $y_i$  $(1\leq i \leq 10)$ are listed in {\em Table 5}.
\end{theorem}

{\bf Proof.} By (2.3), (2.5) and Theorem 2.1, we have $\varphi^{i}(q)\varphi^{j}(q^2)\varphi^{k}(q^3)\varphi^{l}(q^6) \in M_4(\Gamma_0(24),\chi_6)$.
Appealing to Theorem 3.3(c) and using MAPLE,  we obtain the asserted results.
\eop

We now present formulas for $N(1^i, 2^j, 3^k, 6^l;n)$ for each of the octonary quadratic forms $(i,j,k,l)$ listed  in Table 1 in Theorems 3.7--3.9.

\begin{theorem} 
Let $n \in \nn$. Let  $(i,j,k,l)$ be any of the octonary quadratic forms  listed in the first column of Table 1. Then
\beqars
N(1^i, 2^j, 3^k, 6^l;n) &=& x_1\sigma_{(3,\chi_0,\chi_4)}(n) + x_2\sigma_{(3,\chi_4,\chi_0)}(n) + x_3\sigma_{(3,\chi_0,\chi_4)}(n/3) + x_4\sigma_{(3,\chi_4,\chi_0)}(n/3)\\
&& + \sum_{i=1}^{10}y_ia_i(n),
\eeqars
where $x_i$ $(1\leq i \leq 4)$ and $y_i$ $(1\leq i \leq 10)$ are listed in {\em Table 3}.
\end{theorem}

{\bf Proof.} This follows from (2.2), (2.11), (2.18), (3.1) and Theorem 3.4.
\eop

\begin{theorem} 
Let $n \in \nn$. Let  $(i,j,k,l)$ be any of the octonary quadratic forms  listed in the second column of Table 1. Then
\beqars
N(1^i, 2^j, 3^k, 6^l;n) &=& x_1\sigma_{(3,\chi_0,\chi_5)}(n) + x_2\sigma_{(3,\chi_2,\chi_3)}(n) + x_3\sigma_{(3,\chi_3,\chi_2)}(n) + x_4\sigma_{(3,\chi_5,\chi_0)}(n)\\
&& + x_5\sigma_{(3,\chi_0,\chi_5)}(n/2) + x_6\sigma_{(3,\chi_3,\chi_2)}(n/2) + \sum_{i=1}^{8}y_ib_i(n),
\eeqars
where $x_i$ $(1\leq i \leq 6)$ and $y_i$ $(1\leq i \leq 8)$ are listed in {\em Table 4}.
\end{theorem}

{\bf Proof.} The assertions follow from (2.2), (2.12), (2.15), (2.17), (2.19), (3.2) and Theorem 3.5.
\eop

\begin{theorem} 
Let $n \in \nn$. Let  $(i,j,k,l)$ be any of the octonary quadratic forms  listed in the third column of Table 1. Then
\beqars
N(1^i, 2^j, 3^k, 6^l;n) &=& x_1\sigma_{(3,\chi_0,\chi_6)}(n) + x_2\sigma_{(3,\chi_1,\chi_3)}(n) + x_3\sigma_{(3,\chi_3,\chi_1)}(n) + x_4\sigma_{(3,\chi_6,\chi_0)}(n) \\
&&+ \sum_{i=1}^{10}y_ic_i(n),
\eeqars
where $x_i$ $(1\leq i \leq 4)$ and $y_i$ $(1\leq i \leq 10)$ are listed in {\em Table 5}.
\end{theorem}

{\bf Proof.} The assertions follow from (2.2), (2.13), (2.14), (2.16),  (2.20), (3.3) and Theorem 3.6.
\eop

\begin{landscape}
\setlength\LTcapwidth{\textwidth}
 \setlength{\tabcolsep}{3pt}
\setlength\extrarowheight{3pt}
{\footnotesize
\begin{longtable}{ccrrrrrrrrrrrrrr}
\caption{Values of  $x_i$ ($1\leq i\leq 4$) and  $y_j$  ($1\leq j\leq 10$)  for Theorems 3.4 and 3.7} \\
\hline\hline
\mbox{}
\hspace{-5mm}
No. & $(i,j,k,l)$ & $x_1$ & $x_2$ & $x_3$  & $x_4$  & $y_1$ & $y_2$ & $y_3$ & $y_4$ & $y_5$ & $y_6$ & $y_7$ & $y_8$ & $y_9$ & $y_{10}$\\[1mm]
\hline
\endfirsthead
\hline
\hline
No. & $(i,j,k,l)$ & $x_1$ & $x_2$ & $x_3$  & $x_4$  & $y_1$ & $y_2$ & $y_3$ & $y_4$ & $y_5$ & $y_6$ & $y_7$ & $y_8$ & $y_9$ & $y_{10}$\\[1mm]
\hline
\endhead
\hline 
\endfoot
\hline
\endlastfoot
1& $(0,2,1,5)$  & $\frac{4}{451}$  & $\frac{32}{451}$ & $\frac{78}{451}$ &  $-\frac{624}{451}$ & $-\frac{896}{451}$ & $\frac{1544}{451}$  & $\frac{720}{451}$ & $\frac{2480}{451}$ &  $\frac{8752}{451}$&  $\frac{12752}{451}$ & $-\frac{3152}{451}$ & $\frac{5312}{451}$  & $\frac{860}{451}$ & $-\frac{320}{41}$\\[1mm]
2& $(0,2,3,3)$  & $\frac{4}{451}$  & $\frac{64}{451}$ & $\frac{78}{451}$ &  $-\frac{1248}{451}$ & $-\frac{1604}{451}$ & $\frac{1288}{451}$  & $\frac{2316}{451}$ & $\frac{1572}{451}$ &  $\frac{18916}{451}$&  $\frac{444}{11}$ & $-\frac{7164}{451}$ & $\frac{14112}{451}$  & $\frac{1536}{451}$ & $-\frac{6320}{451}$\\[1mm]
3& $(0,2,5,1)$  & $\frac{4}{451}$  & $\frac{128}{451}$ & $\frac{78}{451}$ &  $-\frac{2496}{451}$ & $-\frac{1216}{451}$ & $\frac{776}{451}$  & $\frac{3704}{451}$ & $-\frac{5656}{451}$ &  $\frac{19400}{451}$&  $\frac{27304}{451}$ & $-\frac{35032}{451}$ & $\frac{46144}{451}$  & $\frac{1084}{451}$ & $-\frac{4704}{451}$\\[1mm]
4& $(0,4,1,3)$  & $-\frac{8}{451}$  & $\frac{64}{451}$ & $\frac{90}{451}$ &  $\frac{720}{451}$ & $-\frac{1532}{451}$ & $\frac{3104}{451}$  & $-\frac{1780}{451}$ & $\frac{14532}{451}$ &  $\frac{4948}{451}$&  $\frac{7948}{451}$ & $\frac{29780}{451}$ & $-\frac{37088}{451}$  & $\frac{36}{11}$ & $-\frac{6128}{451}$\\[1mm]
5& $(0,4,3,1)$  & $-\frac{8}{451}$  & $\frac{128}{451}$ & $\frac{90}{451}$ &  $\frac{1440}{451}$ & $-\frac{6392}{451}$ & $\frac{2592}{451}$  & $-\frac{2360}{451}$ & $\frac{30264}{451}$ &  $\frac{27736}{451}$&  $-\frac{11816}{451}$ & $\frac{81944}{451}$ & $-\frac{104768}{451}$  & $\frac{6272}{451}$ & $-\frac{25504}{451}$\\[1mm]
6& $(0,6,1,1)$  & $\frac{28}{451}$  & $\frac{224}{451}$ & $\frac{54}{451}$ &  $-\frac{432}{451}$ & $\frac{9800}{451}$ & $\frac{3592}{451}$  & $-\frac{3816}{451}$ & $-\frac{23640}{451}$ &  $-\frac{58456}{451}$&  $\frac{37112}{451}$ & $-\frac{84056}{451}$ & $\frac{108032}{451}$  & $-\frac{10052}{451}$ & $\frac{39648}{451}$\\[1mm]
7 & $(1,1,0,6)$  & $-\frac{2}{451}$  & $\frac{16}{451}$ & $\frac{84}{451}$ &  $\frac{672}{451}$ & $-\frac{128}{41}$ & $\frac{776}{451}$  & $\frac{580}{451}$ & $\frac{6708}{451}$ &  $\frac{12020}{451}$ &  $\frac{3996}{451}$ & $\frac{18228}{451}$ & $-\frac{17472}{451}$  & $\frac{56}{11}$ & $-\frac{512}{41}$\\[1mm]
8& $(1,1,2,4)$  & $-\frac{2}{451}$  & $\frac{32}{451}$ & $\frac{84}{451}$ &  $\frac{1344}{451}$ & $-\frac{1680}{451}$ & $\frac{648}{451}$  & $\frac{1296}{451}$ & $\frac{9616}{451}$ &  $\frac{16528}{451}$&  $\frac{4016}{451}$ & $\frac{22864}{451}$ & $-\frac{26848}{451}$  & $\frac{232}{41}$ & $-\frac{6704}{451}$\\[1mm]
9& $(1,1,4,2)$  & $-\frac{2}{451}$  & $\frac{64}{451}$ & $\frac{84}{451}$ &  $\frac{2688}{451}$ & $-\frac{2224}{451}$ & $\frac{392}{451}$  & $\frac{84}{41}$ & $\frac{13628}{451}$ &  $\frac{16524}{451}$&  $-\frac{1356}{451}$ & $\frac{37548}{451}$ & $-\frac{45600}{451}$  & $\frac{3064}{451}$ & $-\frac{8848}{451}$\\[1mm]
10 & $(1,1,6,0)$  & $-\frac{2}{451}$  & $\frac{128}{451}$ & $\frac{84}{451}$ &  $\frac{5376}{451}$ & $\frac{296}{451}$ & $-\frac{120}{451}$  & $-\frac{1624}{451}$ & $\frac{12632}{451}$ &  $-\frac{6936}{451}$ &  $-\frac{2248}{41}$ & $\frac{57896}{451}$ & $-\frac{68672}{451}$  & $\frac{480}{451}$ & $-\frac{5920}{451}$\\[1mm]
11 & $(1,3,0,4)$  & $\frac{10}{451}$  & $\frac{80}{451}$ & $\frac{72}{451}$ &  $-\frac{576}{451}$ & $-\frac{2732}{451}$ & $\frac{2056}{451}$  & $\frac{4096}{451}$ & $\frac{9480}{451}$ &  $\frac{34344}{451}$&  $\frac{23352}{451}$ & $\frac{4584}{451}$ & $-\frac{5088}{451}$  & $\frac{3544}{451}$ & $-\frac{10768}{451}$\\[1mm]
12& $(1,3,2,2)$  & $\frac{10}{451}$  & $\frac{160}{451}$ & $\frac{72}{451}$ &  $-\frac{1152}{451}$ & $-\frac{812}{451}$ & $\frac{1416}{451}$  & $\frac{4396}{451}$ & $-\frac{252}{451}$ &  $\frac{2428}{41}$&  $68$ & $-\frac{16844}{451}$ & $\frac{32000}{451}$  & $\frac{1544}{451}$ & $-\frac{3008}{451}$\\[1mm]
13 & $(1,3,4,0)$  & $\frac{10}{451}$  & $\frac{320}{451}$ & $\frac{72}{451}$ &  $-\frac{2304}{451}$ & $\frac{6636}{451}$ & $\frac{136}{451}$  & $\frac{3192}{451}$ & $-\frac{35952}{451}$ &  $-\frac{19232}{451}$&  $\frac{25456}{451}$ & $-\frac{104800}{451}$ & $\frac{135040}{451}$  & $-\frac{6064}{451}$ & $\frac{26944}{451}$\\[1mm]
14& $(1,5,0,2)$  & $-\frac{26}{451}$  & $\frac{208}{451}$ & $\frac{108}{451}$ &  $\frac{864}{451}$ & $-\frac{5512}{451}$ & $\frac{2872}{451}$  & $\frac{1964}{451}$ & $\frac{27180}{451}$ &  $\frac{45068}{451}$&  $\frac{16196}{451}$ & $\frac{68044}{451}$ & $-\frac{67072}{451}$  & $\frac{152}{11}$ & $-\frac{22048}{451}$\\[1mm]
15 & $(1,5,2,0)$  & $-\frac{26}{451}$  & $\frac{416}{451}$ & $\frac{108}{451}$ &  $\frac{1728}{451}$ & $-\frac{9704}{451}$ & $\frac{1208}{451}$  & $-\frac{2504}{451}$ & $\frac{38088}{451}$ &  $\frac{36104}{451}$&  $-\frac{1784}{41}$ & $\frac{89608}{451}$ & $-\frac{136480}{451}$  & $\frac{10216}{451}$ & $-\frac{38608}{451}$\\[1mm]
16 & $(1,7,0,0)$    & $\frac{2}{11}$  &  $\frac{16}{11}$  & $0$  &  $0$  & $\frac{12}{11}$  & $\frac{24}{11}$   & $-\frac{56}{11}$ & $-\frac{96}{11}$  &  $-\frac{32}{11}$&  $\frac{64}{11}$  & $-32$ & $32$  & $-\frac{8}{11}$ & $\frac{80}{11}$\\[1mm]
17& $(2,0,1,5)$  & $\frac{4}{451}$  & $\frac{64}{451}$ & $\frac{78}{451}$ &  $-\frac{1248}{451}$ & $\frac{200}{451}$ & $\frac{1288}{451}$  & $\frac{512}{451}$ & $\frac{3376}{451}$ &  $\frac{13504}{451}$&  $\frac{224}{11}$ & $\frac{1856}{451}$ & $-\frac{320}{451}$  & $\frac{1536}{451}$ & $\frac{896}{451}$\\[1mm]
18& $(2,0,3,3)$  & $\frac{4}{451}$  & $\frac{128}{451}$ & $\frac{78}{451}$ &  $-\frac{2496}{451}$ & $\frac{588}{451}$ & $\frac{776}{451}$  & $\frac{1900}{451}$ & $\frac{3364}{451}$ &  $\frac{13988}{451}$&  $\frac{11068}{451}$ & $\frac{2852}{451}$ & $\frac{2848}{451}$  & $\frac{1084}{451}$ & $\frac{2512}{451}$\\[1mm]
19& $(2,0,5,1)$  & $\frac{4}{451}$  & $\frac{256}{451}$ & $\frac{78}{451}$ &  $-\frac{4992}{451}$ & $\frac{288}{41}$ & $-\frac{248}{451}$  & $\frac{2872}{451}$ & $-\frac{9288}{451}$ &  $-\frac{4888}{451}$&  $\frac{13032}{451}$ & $-\frac{15000}{451}$ & $\frac{576}{11}$  & $-\frac{1624}{451}$ & $\frac{12960}{451}$\\[1mm]
20& $(2,2,1,3)$  & $-\frac{8}{451}$  & $\frac{128}{451}$ & $\frac{90}{451}$ &  $\frac{1440}{451}$ & $-\frac{4588}{451}$ & $\frac{2592}{451}$  & $\frac{3052}{451}$ & $\frac{24852}{451}$ &  $\frac{43972}{451}$&  $\frac{15244}{451}$ & $\frac{54884}{451}$ & $-\frac{61472}{451}$  & $\frac{6272}{451}$ & $-\frac{18288}{451}$\\[1mm]
21& $(2,2,3,1)$  & $-\frac{8}{451}$  & $\frac{256}{451}$ &  $\frac{90}{451}$ & $\frac{2880}{451}$ &  $-\frac{5288}{451}$ & $\frac{1568}{451}$ & $\frac{8}{41}$  & $\frac{29256}{451}$ & $\frac{33624}{451}$ &  $-\frac{4440}{451}$&  $\frac{74424}{451}$ & $-\frac{95808}{451}$ & $\frac{6844}{451}$  & $-\frac{20960}{451}$ \\[1mm]
22& $(2,4,1,1)$  & $\frac{28}{451}$  & $\frac{448}{451}$ & $\frac{54}{451}$ &  $-\frac{864}{451}$ & $-\frac{240}{451}$ & $\frac{1800}{451}$  & $\frac{5224}{451}$ & $-\frac{312}{451}$ &  $\frac{33848}{451}$&  $\frac{1000}{11}$ & $-\frac{18824}{451}$ & $\frac{42368}{451}$  & $\frac{1568}{451}$ & $-\frac{288}{451}$\\[1mm]
23& $(3,1,0,4)$  & $\frac{10}{451}$  & $\frac{160}{451}$ & $\frac{72}{451}$ &  $-\frac{1152}{451}$ & $-\frac{6224}{451}$ & $\frac{5024}{451}$  & $\frac{6200}{451}$ & $\frac{26808}{451}$ &  $\frac{6200}{41}$&  $88$ & $\frac{39080}{451}$ & $-\frac{40160}{451}$  & $\frac{8760}{451}$ & $-\frac{24656}{451}$\\[1mm]
24& $(3,1,2,2)$  & $\frac{10}{451}$  & $\frac{320}{451}$ & $\frac{72}{451}$ &  $-\frac{2304}{451}$ & $-\frac{2384}{451}$ & $\frac{3744}{451}$  & $\frac{4996}{451}$ & $\frac{12756}{451}$ &  $\frac{51124}{451}$&  $\frac{41692}{451}$ & $\frac{8852}{451}$ & $\frac{5152}{451}$  & $\frac{4760}{451}$ & $-\frac{9136}{451}$\\[1mm]
25& $(3,1,4,0)$  & $\frac{10}{451}$  & $\frac{640}{451}$ & $\frac{72}{451}$ &  $-\frac{4608}{451}$ & $\frac{8904}{451}$ & $\frac{1184}{451}$  & $\frac{784}{451}$ & $-\frac{31584}{451}$ &  $-\frac{28128}{451}$&  $\frac{11424}{451}$ & $-\frac{82272}{451}$ & $\frac{2688}{11}$  & $-\frac{6848}{451}$ & $\frac{29120}{451}$\\[1mm]
26& $(3,3,0,2)$  & $-\frac{26}{451}$  & $\frac{416}{451}$ & $\frac{108}{451}$ &  $\frac{1728}{451}$ & $-\frac{7900}{451}$ & $\frac{4816}{451}$  & $\frac{6516}{451}$ & $\frac{43500}{451}$ &  $\frac{84812}{451}$&  $\frac{3628}{41}$ & $\frac{95020}{451}$ & $-\frac{93184}{451}$  & $\frac{10216}{451}$ & $-\frac{31392}{451}$\\[1mm]
27& $(3,3,2,0)$  & $-\frac{26}{451}$  & $\frac{832}{451}$ & $\frac{108}{451}$ &  $\frac{3456}{451}$ & $-\frac{5460}{451}$ & $\frac{1488}{451}$  & $-\frac{384}{41}$ & $\frac{27432}{451}$ &  $\frac{14568}{451}$&  $-\frac{29928}{451}$ & $\frac{71400}{451}$ & $-\frac{116544}{451}$  & $\frac{7360}{451}$ & $-\frac{21216}{451}$\\[1mm]
28 & $(3,5,0,0)$    & $\frac{2}{11}$  &  $\frac{32}{11}$  & $0$  &  $0$  & $-\frac{8}{11}$  & $-\frac{16}{11}$  & $-\frac{32}{11}$ & $-\frac{144}{11}$ &  $\frac{160}{11}$&  $\frac{304}{11}$ & $-\frac{736}{11}$ & $\frac{736}{11}$  & $\frac{40}{11}$ & $\frac{16}{11}$\\[1mm]
29& $(4,0,1,3)$  & $-\frac{8}{451}$  & $\frac{256}{451}$ & $\frac{90}{451}$ &  $\frac{2880}{451}$ & $-\frac{10700}{451}$ & $\frac{8784}{451}$  & $\frac{500}{41}$ & $\frac{45492}{451}$ &  $\frac{107588}{451}$&  $\frac{58700}{451}$ & $\frac{76228}{451}$ & $-\frac{81376}{451}$  & $\frac{14060}{451}$ & $-\frac{42608}{451}$\\[1mm]
30& $(4,0,3,1)$  & $-\frac{8}{451}$  & $\frac{512}{451}$ & $\frac{90}{451}$ &  $\frac{5760}{451}$ & $-\frac{936}{41}$ & $\frac{6736}{451}$  & $-\frac{2232}{451}$ & $\frac{48888}{451}$ &  $\frac{88696}{451}$&  $\frac{10312}{451}$ & $\frac{59384}{451}$ & $-\frac{77888}{451}$  & $\frac{13400}{451}$ & $-\frac{40736}{451}$\\[1mm]
31& $(4,2,1,1)$  & $\frac{28}{451}$  & $\frac{896}{451}$ & $\frac{54}{451}$ &  $-\frac{1728}{451}$ & $-\frac{5888}{451}$ & $\frac{5432}{451}$  & $\frac{8872}{451}$ & $\frac{24696}{451}$ &  $\frac{88568}{451}$&  $\frac{70424}{451}$ & $\frac{25048}{451}$ & $-\frac{2368}{451}$  & $\frac{8572}{451}$ & $-\frac{22432}{451}$\\[1mm]
32& $(5,1,0,2)$  & $-\frac{26}{451}$  & $\frac{832}{451}$ & $\frac{108}{451}$ &  $\frac{3456}{451}$ & $-\frac{21696}{451}$ & $\frac{12312}{451}$  & $\frac{1748}{41}$ & $\frac{97788}{451}$ &  $\frac{207596}{451}$&  $\frac{108980}{451}$ & $\frac{192268}{451}$ & $-\frac{188704}{451}$  & $\frac{25400}{451}$ & $-\frac{86160}{451}$\\[1mm]
33& $(5,1,2,0)$  & $-\frac{26}{451}$  & $\frac{1664}{451}$ & $\frac{108}{451}$ &  $\frac{6912}{451}$ & $-\frac{5992}{451}$ & $\frac{5656}{451}$  & $-\frac{4056}{451}$ & $\frac{27768}{451}$ &  $\frac{29224}{451}$&  $-\frac{14456}{451}$ & $\frac{63848}{451}$ & $-\frac{105536}{451}$  & $\frac{8864}{451}$ & $-\frac{29728}{451}$\\[1mm]
34 & $(5,3,0,0)$  & $\frac{2}{11}$  & $\frac{64}{11}$ & $0$ &  $0$ & $-\frac{4}{11}$ & $-\frac{8}{11}$  & $-\frac{72}{11}$ & $-\frac{240}{11}$ &  $\frac{192}{11}$&  $\frac{432}{11}$ & $-\frac{1152}{11}$ & $\frac{1152}{11}$  & $\frac{48}{11}$ & $\frac{64}{11}$\\[1mm]
35& $(6,0,1,1)$  & $\frac{28}{451}$  & $\frac{1792}{451}$ & $\frac{54}{451}$ &  $-\frac{3456}{451}$ & $-\frac{17184}{451}$ & $\frac{12696}{451}$  & $\frac{30600}{451}$ & $\frac{4824}{41}$ &  $\frac{183576}{451}$&  $\frac{18312}{41}$ & $\frac{170520}{451}$ & $-\frac{3648}{11}$  & $\frac{20776}{451}$ & $-\frac{66720}{451}$\\[1mm]
36 & $(7,1,0,0)$  & $\frac{2}{11}$  & $\frac{128}{11}$ & $0$ &  $0$ & $-\frac{40}{11}$ & $-\frac{80}{11}$  & $\frac{112}{11}$ & $\frac{96}{11}$ &  $\frac{256}{11}$ &  $\frac{160}{11}$ & $\frac{128}{11}$ & $-\frac{128}{11}$  & $\frac{64}{11}$ & $-\frac{192}{11}$
\end{longtable}
}
\end{landscape}

\begin{landscape}
\setlength\LTcapwidth{\textwidth}
\setlength\extrarowheight{3pt}
{\small\begin{center}
\begin{longtable}{c c r r r r r r r r r r r r r r}
\caption{Values of  $x_i$ ($1\leq i\leq 6$) and  $y_j$  ($1\leq j\leq 8$)  for Theorems 3.5 and 3.8} \\
\hline\hline
No. & $(i,j,k,l)$ & $x_1$ & $x_2$ & $x_3$  & $x_4$  & $x_5$  & $x_6$  & $y_1$ & $y_2$ & $y_3$ & $y_4$ & $y_5$ & $y_6$ & $y_7$ & $y_8$\\[1mm]
\hline
\endfirsthead
\hline
\hline
No. & $(i,j,k,l)$ & $x_1$ & $x_2$ & $x_3$  & $x_4$  & $x_5$  & $x_6$  & $y_1$ & $y_2$ & $y_3$ & $y_4$ & $y_5$ & $y_6$ & $y_7$ & $y_8$\\[1mm]
\hline
\endhead
\hline
\endfoot
\hline
\endlastfoot
1 &$(0,1,2,5)$  & $0$  & $-\frac{1}{23}$ & $0$ &  $\frac{1}{23}$ & $\frac{1}{23}$ & $\frac{1}{23}$  & $0$ & $\frac{28}{23}$ & $\frac{64}{23}$ &
 $0$ & $-\frac{64}{23}$&  $\frac{312}{23}$ & $\frac{128}{23}$ & $\frac{56}{23}$   \\[1mm]
2 &$(0,1,4,3)$  & $0$  & $-\frac{2}{23}$ & $0$ &  $\frac{2}{23}$ & $\frac{1}{23}$ & $\frac{1}{23}$  & $0$ & $\frac{12}{23}$ &  $\frac{128}{23}$&  $0$ & $-\frac{128}{23}$ & $\frac{656}{23}$  & $\frac{256}{23}$ & $\frac{144}{23}$\\[1mm]
3 &$(0,1,6,1)$  & $0$  & $-\frac{4}{23}$ & $0$ &  $\frac{4}{23}$ & $\frac{1}{23}$ & $\frac{1}{23}$  & $0$ & $-\frac{20}{23}$ &  $\frac{164}{23}$&  $0$ & $-\frac{440}{23}$ & $\frac{976}{23}$  & $\frac{328}{23}$ & $\frac{320}{23}$\\[1mm]
4 &$(0,3,2,3)$  & $0$  & $\frac{1}{23}$ & $0$ &  $\frac{3}{23}$ & $\frac{1}{23}$ & $-\frac{3}{23}$  & $-\frac{4}{23}$ & $\frac{124}{23}$ &  $\frac{12}{23}$&  $0$ & $\frac{304}{23}$ & $-\frac{152}{23}$  & $-\frac{240}{23}$ & $-\frac{248}{23}$\\[1mm]
5 &$(0,3,4,1)$  & $0$  & $\frac{2}{23}$ & $0$ &  $\frac{6}{23}$ & $\frac{1}{23}$ & $-\frac{3}{23}$  & $-\frac{8}{23}$ & $\frac{108}{23}$ &  $\frac{24}{23}$&  $-\frac{272}{23}$ & $\frac{608}{23}$ & $-\frac{96}{23}$  & $-\frac{480}{23}$ & $-\frac{528}{23}$\\[1mm]
6 &$(0,5,2,1)$  & $0$  & $-\frac{1}{23}$ & $0$ &  $\frac{9}{23}$ & $\frac{1}{23}$ & $\frac{9}{23}$  & $-\frac{8}{23}$ & $\frac{140}{23}$ &  $-\frac{152}{23}$&  $\frac{384}{23}$ & $-\frac{320}{23}$ & $-\frac{744}{23}$  & $\frac{640}{23}$ & $\frac{280}{23}$\\[1mm]
7 &$(1,0,1,6)$  & $0$  & $\frac{1}{23}$ & $0$ &  $\frac{1}{23}$ & $\frac{1}{23}$ & $-\frac{1}{23}$  & $\frac{44}{23}$ & $0$ &  $\frac{20}{23}$&  $8$ & $0$ & $\frac{104}{23}$  & $0$ & $0$\\[1mm]
8 &$(1,0,3,4)$  & $0$  & $\frac{2}{23}$ & $0$ &  $\frac{2}{23}$ & $\frac{1}{23}$ & $-\frac{1}{23}$  & $\frac{42}{23}$ & $0$ &  $\frac{86}{23}$&  $\frac{232}{23}$ & $0$ & $\frac{328}{23}$  & $0$ & $0$\\[1mm]
9 &$(1,0,5,2)$  & $0$  & $\frac{4}{23}$ & $0$ &  $\frac{4}{23}$ & $\frac{1}{23}$ & $-\frac{1}{23}$  & $\frac{38}{23}$ & $0$ &  $\frac{126}{23}$&  $\frac{144}{23}$ & $0$ & $\frac{592}{23}$  & $0$ & $0$\\[1mm]
10 &$(1,0,7,0)$  & $\frac{1}{23}$  &  $\frac{8}{23}$  & $\frac{1}{23}$  &  $\frac{8}{23}$  & $0$  & $0$   & $\frac{28}{23}$ & $0$  &  $\frac{140}{23}$&  $-\frac{224}{23}$  & $0$ & $\frac{672}{23}$  & $0$ & $0$\\[1mm]
11 &$(1,2,1,4)$  & $0$  & $-\frac{1}{23}$ & $0$ &  $\frac{3}{23}$ & $\frac{1}{23}$ & $\frac{3}{23}$  & $\frac{44}{23}$ & $\frac{56}{23}$ &  $\frac{148}{23}$&  $\frac{280}{23}$ & $-\frac{128}{23}$ & $\frac{600}{23}$  & $\frac{256}{23}$ & $\frac{112}{23}$\\[1mm]
12 &$(1,2,3,2)$  & $0$  & $-\frac{2}{23}$ & $0$ &  $\frac{6}{23}$ & $\frac{1}{23}$ & $\frac{3}{23}$  & $\frac{42}{23}$ & $\frac{24}{23}$ &  $\frac{158}{23}$&  $\frac{328}{23}$ & $-\frac{256}{23}$ & $\frac{776}{23}$  & $\frac{512}{23}$ & $\frac{288}{23}$\\[1mm]
13 &$(1,2,5,0)$  & $0$  & $-\frac{4}{23}$ & $0$ &  $\frac{12}{23}$ & $\frac{1}{23}$ & $\frac{3}{23}$  & $\frac{38}{23}$ & $-\frac{40}{23}$ &  $\frac{86}{23}$&  $\frac{240}{23}$ & $-\frac{880}{23}$ & $\frac{208}{23}$  & $\frac{656}{23}$ & $\frac{640}{23}$\\[1mm]
14 &$(1,4,1,2)$  & $0$  & $\frac{1}{23}$ & $0$ &  $\frac{9}{23}$ & $\frac{1}{23}$ & $-\frac{9}{23}$  & $\frac{36}{23}$ & $\frac{128}{23}$ &  $\frac{172}{23}$&  $8$ & $\frac{480}{23}$ & $24$  & $-\frac{224}{23}$ & $-\frac{256}{23}$\\[1mm]
15 &$(1,4,3,0)$  & $0$  & $\frac{2}{23}$ & $0$ &  $\frac{18}{23}$ & $\frac{1}{23}$ & $-\frac{9}{23}$  & $\frac{26}{23}$ & $\frac{64}{23}$ &  $\frac{22}{23}$&  $-\frac{312}{23}$ & $\frac{224}{23}$ & $\frac{104}{23}$  & $-\frac{1184}{23}$ & $-\frac{640}{23}$\\[1mm]
16 &$(1,6,1,0)$  & $0$  & $-\frac{1}{23}$ & $0$ &  $\frac{27}{23}$ & $\frac{1}{23}$ & $\frac{27}{23}$  & $\frac{20}{23}$ & $\frac{24}{23}$ &  $-\frac{132}{23}$&  $-\frac{40}{23}$ & $-\frac{160}{23}$ & $-\frac{360}{23}$  & $\frac{1056}{23}$ & $\frac{48}{23}$\\[1mm]
17 &$(2,1,0,5)$  & $0$  & $\frac{1}{23}$ & $0$ &  $\frac{3}{23}$ & $\frac{1}{23}$ & $-\frac{3}{23}$  & $\frac{88}{23}$ & $\frac{124}{23}$ &  $\frac{104}{23}$&  $16$ & $-\frac{64}{23}$ & $\frac{216}{23}$  & $\frac{128}{23}$ & $-\frac{248}{23}$\\[1mm]
18 &$(2,1,2,3)$  & $0$  & $\frac{2}{23}$ & $0$ &  $\frac{6}{23}$ & $\frac{1}{23}$ & $-\frac{3}{23}$  & $\frac{84}{23}$ & $\frac{108}{23}$ &  $\frac{116}{23}$&  $\frac{464}{23}$ & $\frac{240}{23}$ & $\frac{272}{23}$  & $-\frac{112}{23}$ & $-\frac{160}{23}$\\[1mm]
19 &$(2,1,4,1)$  & $0$  & $\frac{4}{23}$ & $0$ &  $\frac{12}{23}$ & $\frac{1}{23}$ & $-\frac{3}{23}$  & $\frac{76}{23}$ & $\frac{76}{23}$ &  $\frac{48}{23}$&  $\frac{288}{23}$ & $\frac{296}{23}$ & $\frac{16}{23}$  & $-\frac{408}{23}$ & $-\frac{352}{23}$\\[1mm]
20 &$(2,3,0,3)$  & $0$  & $-\frac{1}{23}$ & $0$ &  $\frac{9}{23}$ & $\frac{1}{23}$ & $\frac{9}{23}$  & $\frac{84}{23}$ & $\frac{140}{23}$ &  $\frac{308}{23}$&  $\frac{752}{23}$ & $\frac{48}{23}$ & $\frac{1096}{23}$  & $\frac{272}{23}$ & $\frac{280}{23}$\\[1mm]
21 &$(2,3,2,1)$  & $0$  & $-\frac{2}{23}$ & $0$ &  $\frac{18}{23}$ & $\frac{1}{23}$ & $\frac{9}{23}$  & $\frac{76}{23}$ & $\frac{60}{23}$ &  $\frac{156}{23}$&  $\frac{576}{23}$ & $-\frac{272}{23}$ & $\frac{768}{23}$  & $\frac{912}{23}$ & $\frac{352}{23}$\\[1mm]
22 &$(2,5,0,1)$  & $0$  & $\frac{1}{23}$ & $0$ &  $\frac{27}{23}$ & $\frac{1}{23}$ & $-\frac{27}{23}$  & $\frac{64}{23}$ & $\frac{140}{23}$ &  $\frac{192}{23}$&  $16$ & $\frac{640}{23}$ & $\frac{824}{23}$  & $\frac{192}{23}$ & $-\frac{280}{23}$\\[1mm]
23 &$(3,0,1,4)$  & $0$  & $-\frac{2}{23}$ & $0$ &  $\frac{6}{23}$ & $\frac{1}{23}$ & $\frac{3}{23}$  & $\frac{134}{23}$ & $\frac{208}{23}$ &  $\frac{66}{23}$&  $\frac{696}{23}$ & $\frac{112}{23}$ & $\frac{40}{23}$  & $\frac{144}{23}$ & $-\frac{80}{23}$\\[1mm]
24 &$(3,0,3,2)$  & $0$  & $-\frac{4}{23}$ & $0$ &  $\frac{12}{23}$ & $\frac{1}{23}$ & $\frac{3}{23}$  & $\frac{130}{23}$ & $\frac{144}{23}$ &  $-\frac{6}{23}$&  $\frac{976}{23}$ & $\frac{224}{23}$ & $-\frac{528}{23}$  & $\frac{288}{23}$ & $-\frac{96}{23}$\\[1mm]
25 &$(3,0,5,0)$  & $\frac{1}{23}$  &  $-\frac{8}{23}$  & $-\frac{3}{23}$  & $\frac{24}{23}$ & $0$  &  $0$  & $\frac{124}{23}$ & $\frac{16}{23}$ &  $-\frac{324}{23}$&  $\frac{992}{23}$ & $-\frac{160}{23}$ & $-\frac{2912}{23}$  & $\frac{96}{23}$ & $-\frac{128}{23}$\\[1mm]
26 &$(3,2,1,2)$  & $0$  & $\frac{2}{23}$ & $0$ &  $\frac{18}{23}$ & $\frac{1}{23}$ & $-\frac{9}{23}$  & $\frac{118}{23}$ & $\frac{248}{23}$ &  $\frac{298}{23}$&  $\frac{792}{23}$ & $\frac{592}{23}$ & $\frac{840}{23}$  & $-\frac{80}{23}$ & $-\frac{272}{23}$\\[1mm]
27 &$(3,2,3,0)$  & $0$  & $\frac{4}{23}$ & $0$ &  $\frac{36}{23}$ & $\frac{1}{23}$ & $-\frac{9}{23}$  & $\frac{98}{23}$ & $\frac{120}{23}$ &  $-\frac{94}{23}$&  $-\frac{16}{23}$ & $\frac{80}{23}$ & $-\frac{976}{23}$  & $-\frac{1264}{23}$ & $-\frac{672}{23}$\\[1mm]
28 &$(3,4,1,0)$  & $0$  & $-\frac{2}{23}$ & $0$ &  $\frac{54}{23}$ & $\frac{1}{23}$ & $\frac{27}{23}$  & $\frac{86}{23}$ & $-\frac{16}{23}$ &  $-\frac{126}{23}$&  $\frac{216}{23}$ & $-\frac{688}{23}$ & $\frac{8}{23}$  & $\frac{1008}{23}$ & $\frac{176}{23}$\\[1mm]
29 &$(4,1,0,3)$  & $0$  & $-\frac{2}{23}$ & $0$ &  $\frac{18}{23}$ & $\frac{1}{23}$ & $\frac{9}{23}$  & $\frac{168}{23}$ & $\frac{428}{23}$ &  $\frac{616}{23}$&  $\frac{1312}{23}$ & $\frac{96}{23}$ & $\frac{1872}{23}$  & $\frac{544}{23}$ & $-\frac{16}{23}$\\[1mm]
30 &$(4,1,2,1)$  & $0$  & $-\frac{4}{23}$ & $0$ &  $\frac{36}{23}$ & $\frac{1}{23}$ & $\frac{9}{23}$  & $\frac{152}{23}$ & $\frac{268}{23}$ &  $\frac{220}{23}$&  $\frac{960}{23}$ & $\frac{8}{23}$ & $\frac{848}{23}$  & $\frac{904}{23}$ & $\frac{128}{23}$\\[1mm]
31 &$(4,3,0,1)$  & $0$  & $\frac{2}{23}$ & $0$ &  $\frac{54}{23}$ & $\frac{1}{23}$ & $-\frac{27}{23}$  & $\frac{128}{23}$ & $\frac{300}{23}$ &  $\frac{384}{23}$&  $\frac{1040}{23}$ & $\frac{1280}{23}$ & $\frac{1440}{23}$  & $\frac{384}{23}$ & $-\frac{240}{23}$\\[1mm]
32 &$(5,0,1,2)$  & $0$  & $\frac{4}{23}$ & $0$ &  $\frac{36}{23}$ & $\frac{1}{23}$ & $-\frac{9}{23}$  & $\frac{190}{23}$ & $\frac{672}{23}$ &  $\frac{918}{23}$&  $\frac{720}{23}$ & $\frac{448}{23}$ & $\frac{4176}{23}$  & $\frac{576}{23}$ & $\frac{64}{23}$\\[1mm]
33 &$(5,0,3,0)$  & $\frac{1}{23}$  & $\frac{8}{23}$ & $\frac{9}{23}$ &  $\frac{72}{23}$ & $0$ & $0$  & $\frac{140}{23}$ & $\frac{416}{23}$ &  $\frac{284}{23}$&  $-\frac{1120}{23}$ & $-\frac{832}{23}$ & $\frac{1696}{23}$  & $-\frac{832}{23}$ & $0$\\[1mm]
34 &$(5,2,1,0)$  & $0$  & $-\frac{4}{23}$ & $0$ &  $\frac{108}{23}$ & $\frac{1}{23}$ & $\frac{27}{23}$  & $\frac{126}{23}$ & $\frac{88}{23}$ &  $-\frac{114}{23}$&  $\frac{176}{23}$ & $-\frac{1008}{23}$ & $-\frac{176}{23}$  & $\frac{912}{23}$ & $\frac{64}{23}$\\[1mm]
35 &$(6,1,0,1)$  & $0$  & $\frac{4}{23}$ & $0$ &  $\frac{108}{23}$ & $\frac{1}{23}$ & $-\frac{27}{23}$  & $\frac{164}{23}$ & $\frac{620}{23}$ &  $\frac{1320}{23}$&  $\frac{2016}{23}$ & $\frac{1640}{23}$ & $\frac{4880}{23}$  & $\frac{1320}{23}$ & $-\frac{160}{23}$\\[1mm]
36 &$(7,0,1,0)$  & $\frac{1}{23}$  & $-\frac{8}{23}$ & $-\frac{27}{23}$ &  $\frac{216}{23}$ & $0$ & $0$ & $\frac{140}{23}$ & $\frac{112}{23}$  & $-\frac{84}{23}$ & $\frac{1120}{23}$ &  $-\frac{1120}{23}$&  $-\frac{2912}{23}$ & $\frac{672}{23}$ & $-\frac{896}{23}$ 
\end{longtable}
\end{center}}
\end{landscape}
\noindent

\begin{landscape}
\setlength\LTcapwidth{\textwidth}
\setlength\extrarowheight{3pt}
{\footnotesize\begin{center}
\begin{longtable}{ccrrrrrrrrrrrrrr}
\caption{Values of $x_i$ ($1\leq i\leq 4$) and  $y_j$  ($1\leq j\leq 10$) for Theorems 3.6 and 3.9} \\
\hline\hline
No. & $(i,j,k,l)$ & $x_1$ & $x_2$ & $x_3$  & $x_4$  & $y_1$ & $y_2$ & $y_3$ & $y_4$ & $y_5$ & $y_6$ & $y_7$ & $y_8$ & $y_9$ & $y_{10}$\\[1mm]
\hline
\endfirsthead
\hline
\hline
No. & $(i,j,k,l)$ & $x_1$ & $x_2$ & $x_3$  & $x_4$  & $y_1$ & $y_2$ & $y_3$ & $y_4$ & $y_5$ & $y_6$ & $y_7$ & $y_8$ & $y_9$ & $y_{10}$\\[1mm]
\hline
\endhead
\hline 
\endfoot
\hline
\endlastfoot
1&$(0,1,1,6)$  & $\frac{1}{261}$  & $-\frac{8}{261}$ & $-\frac{1}{261}$ &  $\frac{8}{261}$ & $\frac{392}{261}$ & $\frac{88}{29}$  & $-\frac{148}{87}$ & $-\frac{3136}{261}$ &  $-\frac{136}{29}$&  $\frac{9472}{261}$ & $\frac{2116}{261}$ & $-\frac{3136}{261}$  & $-\frac{392}{261}$ & $-\frac{400}{261}$\\[1mm]
2&$(0,1,3,4)$  & $\frac{1}{261}$  & $-\frac{16}{261}$ & $-\frac{1}{261}$ &  $\frac{16}{261}$ & $\frac{304}{87}$ & $\frac{132}{29}$  & $-\frac{72}{29}$ & $-\frac{2432}{87}$ &  $-\frac{272}{29}$&  $\frac{5600}{87}$ & $\frac{1256}{87}$ & $-\frac{2432}{87}$  & $-\frac{304}{87}$ & $-\frac{308}{87}$\\[1mm]
3&$(0,1,5,2)$  & $\frac{1}{261}$  & $-\frac{32}{261}$ & $-\frac{1}{261}$ &  $\frac{32}{261}$ & $\frac{1952}{261}$ & $\frac{220}{29}$  & $-\frac{700}{87}$ & $-\frac{15616}{261}$ &  $-\frac{776}{29}$&  $\frac{31456}{261}$ & $\frac{10204}{261}$ & $-\frac{15616}{261}$  & $-\frac{1952}{261}$ & $-\frac{68}{9}$\\[1mm]
4 &$(0,1,7,0)$  & $\frac{1}{261}$  & $-\frac{64}{261}$ & $-\frac{1}{261}$ &  $\frac{64}{261}$ & $\frac{448}{29}$ & $\frac{280}{29}$  & $-\frac{672}{29}$ & $-\frac{3584}{29}$ &  $-\frac{2016}{29}$&  $\frac{5824}{29}$ & $\frac{2912}{29}$ & $-\frac{3584}{29}$  & $-\frac{448}{29}$ & $-\frac{336}{29}$\\[1mm]
5&$(0,3,1,4)$  & $\frac{1}{261}$  & $\frac{8}{261}$ & $\frac{1}{87}$ &  $\frac{8}{87}$ & $-\frac{104}{29}$ & $\frac{60}{29}$  & $\frac{176}{29}$ & $32$ &  $\frac{672}{29}$&  $-\frac{224}{29}$ & $-\frac{848}{29}$ & $\frac{1280}{29}$  & $\frac{100}{29}$ & $\frac{100}{29}$\\[1mm]
6&$(0,3,3,2)$  & $\frac{1}{261}$  & $\frac{16}{261}$ & $\frac{1}{87}$ &  $\frac{16}{87}$ & $-\frac{1000}{261}$ & $\frac{128}{87}$  & $\frac{2036}{261}$ & $\frac{736}{29}$ &  $\frac{3368}{87}$&  $-\frac{5536}{261}$ & $-\frac{11500}{261}$ & $\frac{24064}{261}$  & $\frac{932}{261}$ & $\frac{32}{9}$\\[1mm]
7 &$(0,3,5,0)$  & $\frac{1}{261}$  & $\frac{32}{261}$ & $\frac{1}{87}$ &  $\frac{32}{87}$ & $-\frac{376}{87}$ & $\frac{8}{29}$  & $\frac{632}{87}$ & $\frac{352}{29}$ &  $\frac{1328}{29}$&  $-\frac{4192}{87}$ & $-\frac{9544}{87}$ & $\frac{16384}{87}$  & $\frac{332}{87}$ & $\frac{328}{87}$\\[1mm]
8&$(0,5,1,2)$  & $\frac{1}{261}$  & $-\frac{8}{261}$ & $-\frac{1}{29}$ &  $\frac{8}{29}$ & $-\frac{1384}{261}$ & $\frac{224}{87}$  & $\frac{1484}{261}$ & $\frac{5696}{87}$ &  $\frac{152}{87}$&  $\frac{1664}{261}$ & $\frac{1964}{261}$ & $-\frac{15680}{261}$  & $\frac{1328}{261}$ & $\frac{1288}{261}$\\[1mm]
9 &$(0,5,3,0)$  & $\frac{1}{261}$  & $-\frac{16}{261}$ & $-\frac{1}{29}$ &  $\frac{16}{29}$ & $-\frac{1376}{87}$ & $-\frac{292}{29}$  & $\frac{1984}{87}$ & $\frac{3840}{29}$ &  $\frac{544}{29}$&  $-\frac{16736}{87}$ & $-\frac{1664}{87}$ & $-\frac{12160}{87}$  & $\frac{1336}{87}$ & $\frac{1316}{87}$\\[1mm]
10 &$(0,7,1,0)$  & $\frac{1}{261}$  & $\frac{8}{261}$ & $\frac{3}{29}$ &  $\frac{24}{29}$ & $\frac{168}{29}$ & $\frac{420}{29}$  & $-\frac{1064}{29}$ & $-\frac{672}{29}$ &  $-\frac{1680}{29}$&  $\frac{5600}{29}$ & $\frac{2744}{29}$ & $\frac{1792}{29}$  & $-\frac{196}{29}$ & $-\frac{196}{29}$\\[1mm]
11 &$(1,0,0,7)$  & $\frac{1}{261}$  & $\frac{8}{261}$ & $\frac{1}{261}$ &  $\frac{8}{261}$ & $\frac{28}{29}$ & $-\frac{28}{29}$  & $\frac{28}{29}$ & $\frac{56}{29}$ &  $\frac{56}{29}$&  $-\frac{168}{29}$ & $-\frac{84}{29}$ & $0$  & $\frac{28}{29}$ & $\frac{28}{29}$\\[1mm]
12&$(1,0,2,5)$  & $\frac{1}{261}$  & $\frac{16}{261}$ & $\frac{1}{261}$ &  $\frac{16}{261}$ & $\frac{1004}{261}$ & $\frac{172}{87}$  & $-\frac{16}{9}$ & $-\frac{4544}{261}$ &  $-\frac{344}{87}$&  $\frac{2560}{87}$ & $\frac{568}{87}$ & $0$  & $-\frac{172}{87}$ & $-\frac{172}{87}$\\[1mm]
13&$(1,0,4,3)$  & $\frac{1}{261}$  & $\frac{32}{261}$ & $\frac{1}{261}$ &  $\frac{32}{261}$ & $\frac{488}{87}$ & $\frac{112}{29}$  & $-\frac{284}{87}$ & $-\frac{2792}{87}$ &  $-\frac{224}{29}$&  $\frac{1736}{29}$ & $\frac{388}{29}$ & $0$  & $-\frac{112}{29}$ & $-\frac{112}{29}$\\[1mm]
14&$(1,0,6,1)$  & $\frac{1}{261}$  & $\frac{64}{261}$ & $\frac{1}{261}$ &  $\frac{64}{261}$ & $\frac{2384}{261}$ & $\frac{664}{87}$  & $-\frac{2672}{261}$ & $-\frac{18128}{261}$ &  $-\frac{1328}{87}$&  $\frac{9808}{87}$ & $\frac{1312}{87}$ & $0$  & $-\frac{664}{87}$ & $-\frac{664}{87}$\\[1mm]
15 &$(1,2,0,5)$  & $\frac{1}{261}$  & $-\frac{8}{261}$ & $-\frac{1}{87}$ &  $\frac{8}{87}$ & $\frac{2036}{261}$ & $\frac{776}{87}$  & $-\frac{1528}{261}$ & $-\frac{4232}{87}$ &  $-\frac{1312}{87}$&  $\frac{26312}{261}$ & $\frac{6776}{261}$ & $-\frac{6272}{261}$  & $-\frac{1528}{261}$ & $-\frac{1544}{261}$\\[1mm]
16&$(1,2,2,3)$  & $\frac{1}{261}$  & $-\frac{16}{261}$ & $-\frac{1}{87}$ &  $\frac{16}{87}$ & $\frac{20}{3}$ & $\frac{200}{29}$  & $-\frac{188}{87}$ & $-\frac{1040}{29}$ &  $-\frac{416}{29}$&  $\frac{8368}{87}$ & $\frac{2092}{87}$ & $-\frac{4864}{87}$  & $-\frac{416}{87}$ & $-\frac{424}{87}$\\[1mm]
17&$(1,2,4,1)$  & $\frac{1}{261}$  & $-\frac{32}{261}$ & $-\frac{1}{87}$ &  $\frac{32}{87}$ & $\frac{2192}{261}$ & $\frac{596}{87}$  & $-\frac{1768}{261}$ & $-\frac{4376}{87}$ &  $-\frac{3208}{87}$&  $\frac{28952}{261}$ & $\frac{14672}{261}$ & $-\frac{31232}{261}$  & $-\frac{1732}{261}$ & $-\frac{1772}{261}$\\[1mm]
18&$(1,4,0,3)$  & $\frac{1}{261}$  & $\frac{8}{261}$ & $\frac{1}{29}$ &  $\frac{8}{29}$ & $-\frac{36}{29}$ & $\frac{92}{29}$  & $\frac{388}{29}$ & $\frac{760}{29}$ &  $\frac{1128}{29}$&  $\frac{536}{29}$ & $-\frac{1516}{29}$ & $\frac{1408}{29}$  & $\frac{84}{29}$ & $\frac{84}{29}$\\[1mm]
19&$(1,4,2,1)$  & $\frac{1}{261}$  & $\frac{16}{261}$ & $\frac{1}{29}$ &  $\frac{16}{29}$ & $\frac{1340}{261}$ & $\frac{692}{87}$  & $-\frac{904}{261}$ & $-\frac{3328}{87}$ &  $\frac{1976}{87}$&  $\frac{21632}{261}$ & $-\frac{3088}{261}$ & $\frac{29440}{261}$  & $-\frac{988}{261}$ & $-\frac{1004}{261}$\\[1mm]
20&$(1,6,0,1)$  & $\frac{1}{261}$  & $-\frac{8}{261}$ & $-\frac{3}{29}$ &  $\frac{24}{29}$ & $\frac{884}{261}$ & $\frac{656}{87}$  & $\frac{80}{261}$ & $-\frac{712}{87}$ &  $\frac{368}{87}$&  $\frac{19592}{261}$ & $\frac{10496}{261}$ & $-\frac{10496}{261}$  & $-\frac{544}{261}$ & $-\frac{656}{261}$\\[1mm]
21&$(2,1,1,4)$  & $\frac{1}{261}$  & $\frac{16}{261}$ & $\frac{1}{87}$ &  $\frac{16}{87}$ & $\frac{1088}{261}$ & $\frac{476}{87}$  & $\frac{992}{261}$ & $-\frac{192}{29}$ &  $\frac{1280}{87}$&  $\frac{11168}{261}$ & $-\frac{4192}{261}$ & $\frac{7360}{261}$  & $-\frac{112}{261}$ & $-\frac{4}{9}$\\[1mm]
22&$(2,1,3,2)$  & $\frac{1}{261}$  & $\frac{32}{261}$ & $\frac{1}{87}$ &  $\frac{32}{87}$ & $\frac{320}{87}$ & $\frac{124}{29}$  & $\frac{284}{87}$ & $-\frac{112}{29}$ &  $\frac{632}{29}$&  $\frac{2768}{87}$ & $-\frac{1540}{87}$ & $\frac{5248}{87}$  & $-\frac{16}{87}$ & $-\frac{20}{87}$\\[1mm]
23&$(2,1,5,0)$  & $\frac{1}{261}$  & $\frac{64}{261}$ & $\frac{1}{87}$ &  $\frac{64}{87}$ & $\frac{704}{261}$ & $-\frac{184}{87}$  & $-\frac{472}{261}$ & $-\frac{416}{29}$ &  $\frac{1040}{87}$&  $-\frac{9952}{261}$ & $-\frac{14872}{261}$ & $\frac{32512}{261}$  & $\frac{80}{261}$ & $\frac{1096}{261}$\\[1mm]
24&$(2,3,1,2)$  & $\frac{1}{261}$  & $-\frac{16}{261}$ & $-\frac{1}{29}$ &  $\frac{16}{29}$ & $\frac{712}{87}$ & $\frac{288}{29}$  & $\frac{244}{87}$ & $-\frac{800}{29}$ &  $-\frac{152}{29}$&  $\frac{11104}{87}$ & $\frac{2164}{87}$ & $-\frac{6592}{87}$  & $-\frac{404}{87}$ & $-\frac{424}{87}$\\[1mm]
25&$(2,3,3,0)$  & $\frac{1}{261}$  & $-\frac{32}{261}$ & $-\frac{1}{29}$ &  $\frac{32}{29}$ & $\frac{824}{261}$ & $\frac{56}{87}$  & $\frac{272}{261}$ & $-\frac{496}{87}$ &  $-\frac{3760}{87}$&  $\frac{560}{261}$ & $\frac{4064}{261}$ & $-\frac{44672}{261}$  & $-\frac{28}{261}$ & $-\frac{128}{261}$\\[1mm]
26&$(2,5,1,0)$  & $\frac{1}{261}$  & $\frac{16}{261}$ & $\frac{3}{29}$ &  $\frac{48}{29}$ & $\frac{2096}{261}$ & $\frac{644}{87}$  & $-\frac{4504}{261}$ & $-\frac{6016}{87}$ &  $-\frac{112}{87}$&  $\frac{19616}{261}$ & $\frac{10664}{261}$ & $\frac{28864}{261}$  & $-\frac{1528}{261}$ & $-\frac{1580}{261}$\\[1mm]
27 &$(3,0,0,5)$  & $\frac{1}{261}$  & $-\frac{16}{261}$ & $-\frac{1}{87}$ &  $\frac{16}{87}$ & $\frac{20}{3}$ & $\frac{316}{29}$  & $\frac{160}{87}$ & $-\frac{576}{29}$ &  $\frac{280}{29}$&  $\frac{6976}{87}$ & $-\frac{344}{87}$ & $\frac{704}{87}$  & $-\frac{68}{87}$ & $-\frac{76}{87}$\\[1mm]
28&$(3,0,2,3)$  & $\frac{1}{261}$  & $-\frac{32}{261}$ & $-\frac{1}{87}$ &  $\frac{32}{87}$ & $\frac{104}{261}$ & $\frac{248}{87}$  & $\frac{3452}{261}$ & $\frac{3976}{87}$ &  $\frac{3056}{87}$&  $-\frac{4456}{261}$ & $-\frac{9340}{261}$ & $\frac{2176}{261}$  & $\frac{1400}{261}$ & $\frac{1360}{261}$\\[1mm]
29&$(3,0,4,1)$  & $\frac{1}{261}$  & $-\frac{64}{261}$ & $-\frac{1}{87}$ &  $\frac{64}{87}$ & $-\frac{352}{29}$ & $-\frac{384}{29}$  & $32$ & $\frac{4432}{29}$ &  $\frac{1104}{29}$&  $-\frac{6832}{29}$ & $-\frac{1840}{29}$ & $\frac{256}{29}$  & $\frac{512}{29}$ & $\frac{504}{29}$\\[1mm]
30&$(3,2,0,3)$  & $\frac{1}{261}$  & $\frac{16}{261}$ & $\frac{1}{29}$ &  $\frac{16}{29}$ & $\frac{1340}{261}$ & $\frac{1040}{87}$  & $\frac{4316}{261}$ & $\frac{848}{87}$ &  $\frac{4064}{87}$&  $\frac{25808}{261}$ & $-\frac{14572}{261}$ & $\frac{12736}{261}$  & $\frac{56}{261}$ & $\frac{40}{261}$\\[1mm]
31&$(3,2,2,1)$  & $\frac{1}{261}$  & $\frac{32}{261}$ & $\frac{1}{29}$ &  $\frac{32}{29}$ & $\frac{512}{87}$ & $\frac{276}{29}$  & $\frac{248}{87}$ & $-\frac{440}{29}$ &  $\frac{1112}{29}$&  $\frac{6584}{87}$ & $-\frac{1648}{87}$ & $\frac{9856}{87}$  & $-\frac{100}{87}$ & $-\frac{4}{3}$\\[1mm]
32&$(3,4,0,1)$  & $\frac{1}{261}$  & $-\frac{16}{261}$ & $-\frac{3}{29}$ &  $\frac{48}{29}$ & $\frac{1108}{87}$ & $\frac{436}{29}$  & $-\frac{200}{87}$ & $-\frac{1472}{29}$ &  $-\frac{56}{29}$&  $\frac{17920}{87}$ & $\frac{6208}{87}$ & $-\frac{6208}{87}$  & $-\frac{716}{87}$ & $-\frac{772}{87}$\\[1mm]
33&$(4,1,1,2)$  & $\frac{1}{261}$  & $-\frac{32}{261}$ & $-\frac{1}{29}$ &  $\frac{32}{29}$ & $\frac{2912}{261}$ & $\frac{1796}{87}$  & $\frac{3404}{261}$ & $-\frac{1888}{87}$ &  $\frac{2504}{87}$&  $\frac{54848}{261}$ & $-\frac{1156}{261}$ & $-\frac{11264}{261}$  & $-\frac{1072}{261}$ & $-\frac{1172}{261}$\\[1mm]
34&$(4,1,3,0)$  & $\frac{1}{261}$  & $-\frac{64}{261}$ & $-\frac{1}{29}$ &  $\frac{64}{29}$ & $\frac{32}{29}$ & $-\frac{56}{29}$  & $\frac{160}{29}$ & $\frac{640}{29}$ &  $-\frac{672}{29}$&  $-\frac{2112}{29}$ & $-\frac{1248}{29}$ & $-\frac{3072}{29}$  & $\frac{144}{29}$ & $\frac{240}{29}$\\[1mm]
35&$(4,3,1,0)$  & $\frac{1}{261}$  & $\frac{32}{261}$ & $\frac{3}{29}$ &  $\frac{96}{29}$ & $\frac{392}{87}$ & $\frac{152}{29}$  & $-\frac{904}{87}$ & $-\frac{960}{29}$ &  $16$&  $-\frac{64}{87}$ & $-\frac{8}{3}$ & $\frac{12544}{87}$  & $-\frac{4}{87}$ & $-\frac{56}{87}$\\[1mm]
36 &$(5,0,0,3)$  & $\frac{1}{261}$  & $\frac{32}{261}$ & $\frac{1}{29}$ &  $\frac{32}{29}$ & $\frac{2600}{87}$ & $\frac{1552}{29}$  & $\frac{596}{87}$ & $-\frac{5080}{29}$ &  $\frac{416}{29}$&  $\frac{51128}{87}$ & $\frac{92}{87}$ & $-\frac{1280}{87}$  & $-\frac{1840}{87}$ & $-\frac{64}{3}$\\[1mm]
37&$(5,0,2,1)$  & $\frac{1}{261}$  & $\frac{64}{261}$ & $\frac{1}{29}$ &  $\frac{64}{29}$ & $\frac{12368}{261}$ & $\frac{5624}{87}$  & $-\frac{14752}{261}$ & $-\frac{31408}{87}$ &  $-\frac{208}{87}$&  $\frac{206000}{261}$ & $\frac{24752}{261}$ & $-\frac{3584}{261}$  & $-\frac{10408}{261}$ & $-\frac{10520}{261}$\\[1mm]
38&$(5,2,0,1)$  & $\frac{1}{261}$  & $-\frac{32}{261}$ & $-\frac{3}{29}$ &  $\frac{96}{29}$ & $\frac{5072}{261}$ & $\frac{2612}{87}$  & $\frac{2216}{261}$ & $-\frac{5560}{87}$ &  $\frac{2936}{87}$&  $\frac{99128}{261}$ & $\frac{18176}{261}$ & $-\frac{18176}{261}$  & $-\frac{3268}{261}$ & $-\frac{3548}{261}$\\[1mm]
39&$(6,1,1,0)$  & $\frac{1}{261}$  & $\frac{64}{261}$ & $\frac{3}{29}$ &  $\frac{192}{29}$ & $\frac{5600}{261}$ & $\frac{2168}{87}$  & $-\frac{7480}{261}$ & $-\frac{16096}{87}$ &  $-\frac{3952}{87}$&  $\frac{68768}{261}$ & $\frac{9992}{261}$ & $\frac{21760}{261}$  & $-\frac{4288}{261}$ & $-\frac{3608}{261}$\\[1mm]
40 &$(7,0,0,1)$  & $\frac{1}{261}$  & $-\frac{64}{261}$ & $-\frac{3}{29}$ &  $\frac{192}{29}$ & $-\frac{672}{29}$ & $0$ & $\frac{4816}{29}$  & $\frac{11536}{29}$ & $\frac{2352}{29}$ &  $-\frac{2800}{29}$&  $-\frac{1792}{29}$ & $\frac{1792}{29}$ & $\frac{896}{29}$  & $\frac{840}{29}$
\end{longtable}
\end{center}}
\end{landscape}

\textbf{\large Acknowledgments.} The research was carried out during the first author's visit to Recep Tayyip  Erdo\u{g}an University, Rize, Turkey.
She would like to thank  Recep Tayyip Erdo\u{g}an University for the hospitality during her stay. 
She thanks the Scientific and Technological Research Council of Turkey (T\"{U}B\.{I}TAK) for their partial financial support. 
The research of the first author was also partially supported by a Discovery Grant
from the Natural Sciences and Engineering Research Council of Canada (RGPIN-418029-2013).

\noindent
Ay\c{s}e Alaca\\
School of Mathematics and Statistics \\
Carleton University, Ottawa \\
Ontario, K1S 5B6, Canada \\
aalaca@math.carleton.ca

\vspace{2mm}

\noindent
M. Nesibe Kesicio\u{g}lu\\
Department of Mathematics \\
Recep Tayyip Erdo\u{g}an University \\
Rize, 53100, Turkey \\
m.nesibe@gmail.com


\begin{thebibliography}{0}

\bibitem{AAW4}
{A. Alaca, \c S. Alaca and K. S. Williams,}
\newblock {\it Fourteen octonary quadratic forms,}
\newblock Int. J. Number Theory {\bf 6} (2010), 37-50.

\bibitem{AK}
{\c S. Alaca and Y. Kesicio\u{g}lu,}
\newblock {\it Representations by certain octonary quadratic forms whose coefficients are $1$, $2$, $3$ and $6$,}
\newblock Int. J. Number Theory {\bf 10} (2014), 133-150.

\bibitem{AK2}
{\c S. Alaca and Y. Kesicio\u{g}lu,}
\newblock {\it Representations by Certain Octonary Quadratic forms with Coefficients $1$, $2$, $3$ and $6$,}
\newblock Integers, {\bf 15} (2015), 1-9.

\bibitem{SW}
{\c S. Alaca and K. S. Williams,}
\newblock {\it The number of representations of a positive integer by certain octonary quadratic forms,}
\newblock Functiones et Approximatio {\bf 43} (2010), 45-54.

\bibitem{GordonSinor}
{B. Gordon and D. Sinor,}
\newblock {\it Multiplicative properties of $\eta$-products,}
\newblock Lecture Notes in Math, vol.1395 Springer-Verlag, New York (1989), 173-200.

\bibitem{jacobi}
{C. G. J. Jacobi,}
\newblock {\it Fundamenta nova theoriae functionum ellipticarum,}
\newblock Borntr\"{a}ger, Regiomonti, 1829.
(Gesammelte Werke, Erster Band, Chelsea Publishing Co., New York, 1969, pp. 49-239.)

\bibitem{Kilford}
{L. J. P. Kilford,}
\newblock {\it Modular Forms, A classical and computational introduction,}
\newblock Imperial College Press, London, 2008.

\bibitem{Kohler}
{G. K\"{o}hler,}
\newblock {\it Eta Products and Theta Series Identities,}
\newblock Springer, 2011.

\bibitem{kokluce-1}
{B. K\"{o}kl\"{u}ce,}
\newblock {\it The representation numbers of three octonary quadratic forms,}
\newblock Int. J. Number Theory {\bf 9} (2013), 505-516.

\bibitem{Ligozat}
{G. Ligozat,}
\newblock {\it Courbes modulaires de genre 1,}
\newblock Bull. Soc. Math. France {\bf 43} (1975), 5-80.

\bibitem{ramak}
{B. Ramakrishnan and B. Sahu,}
\newblock {\it  Evaluation of the convolution sums $ \sum_{l+15m=n}\sigma(l)\sigma(m)$ and 
$ \sum_{3l+5m=n}\sigma(l)\sigma(m)$  and an application,}
\newblock Int. J. Number Theory {\bf 9} (2013), 799-809. 

\bibitem{Serre}
{J. P. Serre,}
\newblock {\it Modular forms of weight one and Galois representations,}
\newblock Algebraic number fields: L -functions and Galois properties (Proc. Sympos., Univ. Durham, Durham, 1975), Academic Press, London, 1977, 193-268.

\bibitem{stein}
{W. Stein,}
\newblock {\it Modular Forms, A Computational Approach,}
\newblock Ameican Mathematical Society, 2007.

\end{thebibliography}
\end{document}